\documentclass[12pt,a4paper]{article}

\newcommand{\sub}{\supseteq}

\newcommand{\tohi}{\varnothing}

\def\ifif {if and only if\ \ }
\def\sub {\subseteq}

\def\set {\setminus }

\def\RE {{\Bbb R}}

\def\N {{\Bbb N}}

\usepackage{amsmath}
\usepackage{amsfonts}
\usepackage{amssymb}
\usepackage{makeidx}
\usepackage{amsthm}
\newcommand{\mathscr}{\mathbf}

\begin{document}
\title{Some generalizations and unifications of $C_{K}(X)$, $C_{\psi}(X)$ and $C_{\infty}(X)$}
\author{A. Taherifar \\Department of Mathematics, Yasouj University, Yasouj , Iran}
\date{ataherifar@mail.yu.ac.ir \[ \textbf{\text{ Dedicated to Professor Azarpanah}}\]}

\maketitle

\theoremstyle{definition}\newtheorem{thm}{Theorem}[section]
\theoremstyle{definition}\newtheorem{cor}[thm]{Corollary}
\theoremstyle{definition}\newtheorem{lem}[thm]{Lemma}
\theoremstyle{definition}\newtheorem{pro}[thm]{Proposition}
\theoremstyle{definition}\newtheorem{defn}[thm]{Definition}
\theoremstyle{definition}\newtheorem{Rem}[thm]{Remark}
\theoremstyle{remark}\newtheorem{rem}[thm]{Remark}
\theoremstyle{definition}\newtheorem{exa}[thm]{Example}
\theoremstyle{definition}\newtheorem{qu}[thm]{question}
{\noindent\bf Abstract.} Let $\cal{P}$ be an open  filter base for
a filter $\cal{F}$ on $X$. We denote by $C^{\cal{P}}(X)$
($C_{\infty\cal{P}}(X)$) the set of all functions $f\in C(X)$
where $Z(f)$ ($\{x: |f(x)|< \frac{1}{n}\})$ contains an element of
$\cal{P}$. First, we observe that every proper subrings in the
sense of Acharyya and Ghosh (Topology Proc. 2010) has such form
and vice versa. After wards, we generalize some well known
theorems about $C_{K}(X), C_{\psi}(X)$ and $C_{\infty}(X)$ for
$C^{\cal{P}}(X)$ and $C_{\infty\cal{P}}(X)$. We observe that
$C_{\infty\cal{P}}(X)$ may not be an ideal of $C(X)$. It is  shown
that $C_{\infty\cal{P}}(X)$ is an ideal of $C(X)$ and for each
$F\in\cal{F}$, $X\setminus \overline{F}$ is bounded \ifif the set
of non-cluster points of the filter $\cal{F}$ is bounded. By this
result, we investigate topological spaces for which
$C_{\infty\cal{P}}(X)$ is an ideal of $C(X)$ whenever
$\cal{P}$=$\{A\subsetneq X$: $A$ is open and $X\setminus A$ is
bounded $\}$ (resp., $\cal{P}$=$\{A\subsetneq X$: $X\setminus A$
is finite $\}$). Moreover, we prove that $C^{\cal{P}}(X)$ is an
essential (resp., free) ideal \ifif the set $\{V:$ $V$ is open and
$X\setminus V\in\mathcal{F}\}$ is a $\pi$-base for $X$ (resp.,
$\mathcal{F}$ has no cluster point). Finally, the filter $\cal{F}$
for which $C_{\infty\cal{P}}(X)$ is a regular ring (resp.,
$z$-ideal) is characterized.

\vspace{0.3 cm}\noindent{\it Keywords:} local space, bounded subset, $z$-ideal, regular ring, essential ideal, $\cal{F}$-$CG_{\delta}$ subset.

\vspace{0.3 cm}\noindent{\it 2010 Mathematical Subject
Classification: $54C40$.}
\section{Introduction}
In this paper, $X$ assumed to be a completely regular Hausdorff
space. $C(X)$($C^{*}(X))$ stands for the ring of all real valued
(bounded) continuous functions on $X$. A subcollection
$\mathcal{B}$ of a filter $\mathcal{F}$ is a $\emph{filter}$
$\emph{base}$ for $\mathcal{F}$ \ifif each element of
$\mathcal{F}$ contains some element of $\mathcal{B}$. A nonempty
collection $\cal{B}$ of nonempty subsets of the space $X$ is a
filter base for some filter on $X$ \ifif the intersection of any
two elements of $\cal{B}$ contains an element of $\cal{B}$. If any
element of $\cal{B}$ is an open subset, we say $\cal{B}$ is an
open  filter base. We shall say that a point $p\in X$ is a cluster
point of $\cal{F}$ if every neighborhood of $p$ meets every member
of $\cal{F}$, in other words,
$p\in\bigcap_{F\in\cal{F}}\overline{F}$. Kohls proved in \cite{K}
 that the intersection of all free maximal ideals in
$C^{*}(X)$ is precisely the set $C_{\infty}(X)$ consists of all
continuous functions $f\in C(X)$ which vanish at infinity, in the
sense that $\{x: |f(x)|\geq\frac{1}{n}\}$ is compact for each
$n\in\Bbb N$. Kohls  also shown in \cite{K} that the $C_{K}(X)$
consists of all functions in $C(X)$ with compact support is the
intersection of all free ideals in $C(X)$ and all free ideals in
$C^{*}(X)$. It is well known that $C_{K}(X)$ is an ideal of $C(X)$
and it is easy to see that  $C_{\infty}(X)$ is an ideal of
$C^{*}(X)$ but not in $C(X)$, see \cite{AS}. An ideal $I$ of
$C(X)$ is called essential ideal if $I\cap (f)\neq 0$, for every
non-zero element $f\in C(X)$.  An ideal $I$ of $C(X)$ is called
$z$-ideal if $Z(f)\sub Z(g)$, $f\in I$, then $g\in I$. In
particular, maximal ideals, minimal prime ideals, and most of
familiar ideals in $C(X)$ are $z$-ideals, see \cite{ART} and
\cite{AT}. \\\indent Our aim of this paper is to reveal some
important properties of a special kind of generalized form of
$C_{K}(X)$ and $C_{\infty}(X))$, which is denoted by
$C^{\cal{P}}(X)$ and $C_{\infty\cal{P}}(X)$. In section $2$, some
examples of these subrings are given and  we prove that for any
open filter base $\mathcal{P}$, there is an ideal
$\mathcal{P}^{\prime}$ of closed sets such that
$C^{\cal{P}}(X)=C_{\cal{P}^{\prime}}(X)$ and
$C_{\infty\cal{P}}(X)=C_{\infty}^{\cal{P}^{\prime}}(X)$ and vice
versa whenever  $C_{\infty}^{\cal{P}^{\prime}}(X)$ is a proper
subsring pf $C(X)$ (see \cite{AK}). It is shown that
$C^{\cal{P}}(X)$ is a free ideal \ifif $\cal{F}$ has no cluster
point. Consequently, we observe that $X$ is a local space (i,e.,
there is an open filter base $\cal{P}$ for a filter $\cal{F}$
which $\cal{F}$ has no cluster point) \ifif there is an open
filter base $\cal{P}$ such that $C^{\cal{P}}(X)$ is a free ideal.
In this section, we show that $C^{\cal{P}}(X)$ (resp.,
$C_{\infty\cal{P}}(X)$) is a zero ideal \ifif any element of
$\cal{P}$ ($\cal{F}$) is dense in $X$. A subset $A$ of $X$ is
called an $\cal{F}$-$CG_{\delta}$, if
$A=\bigcap_{i=1}^{\infty}A_{i}$, where each $A_{i}$ is an open
subset, $X\set A_{i}$ and $\overline{A_{i+1}}$ are completely
separated and each $A_{i}\in \cal{F}$. We prove that that
$C^{\cal{P}}(X)=C_{\infty\cal{P}}(X)$ \ifif every closed
$\cal{F}$-$CG_{\delta}$ is an element of $\cal{F}$. We give an
example of an open filter base $\cal{P}$ for a filter
$\mathcal{F}$ such that  $C_{\infty\cal{P}}(X)$ is not an ideal of
$C(X)$. It is also shown that $C_{\infty\cal{P}}(X)$ is an ideal
of $C(X)$ and for any $F\in\cal{F}$, $X\set \overline{F}$ is
bounded \ifif the set of non-cluster points of the filter
$\cal{F}$ is bounded which is a generalization of Theorem 1.3 in
\cite{AS}. Consequently, if $X$ is a pseudocompact space, then for
any open filter base $\cal{P}$, $C_{\infty\cal{P}}(X)$ is an ideal
of $C(X)$. \\\indent Section 3 is devoted to the essentiality of
$C^{\cal{P}}(X)$ and ideals in $C_{\infty\cal{P}}(X)$. In
\cite{A}, it was proved that an ideal $I$ of $C(X)$ is an
essential ideal \ifif $\bigcap Z[I]$ does not contain an open
subset (i.e., int$\bigcap Z[I]=\tohi$).
 Azarpanah in \cite{AZ}, proved that $C_{K}(X)$ is an essential ideal
\ifif $X$ is an almost locally compact (i.e, every non-empty open
set of $X$ contains a non-empty open set with compact closure). We
generalize these results for $C^{\cal{P}}(X)$ and
$C_{\infty\cal{P}}(X)$. It is proved that an ideal $E$ in
$C_{\infty\cal{P}}(X)$ is an essential ideal \ifif $\bigcap Z[E]$
does not contain a subset $V$ where $X\setminus V\in\mathcal{F}$
and $intV\neq\emptyset$. We also prove that $C^{\cal{P}}(X)$ is an
essential ideal \ifif the set $\{V:$ $V$ is open and $X\setminus
V\in\mathcal{F}\}$ is a $\pi$-base for $X$. \\\indent In section
$4$, it is shown that $C_{\infty\cal{P}}(X)$ is a $z$-ideal \ifif
every cozero-set containing  a closed $\cal{F}$-$CG_{\delta}$ is
an element of $\cal{F}$. Also we show that $C_{\infty\cal{P}}(X)$
is a regular ring (in the sense of Von Neumann) \ifif every closed
$\cal{F}$-$CG_{\delta}$ is an open subset and belongs to
$\cal{F}$. Finally, we see that if $\cal{P}$=\{$A$: $A$ is open
and $X\setminus A$ is Lindel\"{o}f$\}$ and $X$ is a
non-Lindel\"{o}f space, then $\cal{P}$ is an open filter and
$C_{\infty\cal{P}}(X)$ is a regular ring \ifif every closed
$\cal{P}$-$CG_{\delta}$ is an open subset.

\section{ $C^{\cal{P}}(X)$ and $C_{\infty\cal{P}}(X)$}

Let $\cal{P}$ be an open filter base for a filter $\mathcal{F}$ on
topological space $X$, we denote by $C^{\cal{P}}(X)$  the set of
all functions $f$ in $C(X)$ for which $ Z(f)$ contains an element
of
   $\cal{P}$.
  Also $C_{\infty\cal{P}}(X)$ denotes the family of all
functions $f\in C(X)$ for which the set $\{x: |f(x)|<
\frac{1}{n}\}$ contains an element of
 $\cal{P}$, for each $n\in \Bbb N$.
 \\\indent Recall that, for a subset $A$ of $X$, $O_{A}=\{f: A\subseteq int Z(f)\}$.
\begin{lem}\label{2.2}
The following statements hold.
\begin{enumerate}
 \item $C^{\cal{P}}(X)$ is a $z$-ideal
of $C(X)$ contained in $C_{\infty\cal{P}}(X)$.

\item
$C^{\cal{P}}(X)=\sum_{A\in\cal{P}}O_{A}=\bigcup_{A\in\cal{P}}O_{A}$

\item  $C_{\infty\cal{P}}(X)$ is a proper subring of $C(X)$.
\end{enumerate}
\end{lem}
{\noindent\bf Proof.} (1) By definition of $ C^{\cal{P}}(X)$ and
since $\cal{P}$ is a base filter, it is easy to see that
$C^{\cal{P}}(X)$ is  a $z$-ideal and contained in
$C_{\infty\cal{P}}(X)$.

(2) Let $f\in  C^{\cal{P}}(X)$. Then there exists $A\in\cal{P}$
such that $A\subseteq Z(f)$. Hence $A\subseteq intZ(f)$, i.e.,
$f\in O_{A}\subseteq \sum_{A\in\cal{P}}O_{A}$. If $f\in
\sum_{A\in\cal{P}}O_{A}$, then there are $1\leq i\leq n$ and
$f_{i}\in O_{A_{i}}$ such that $f=f_{1}+...+f_{n}$, thus
$\bigcap_{i=1}^{n}A_{i}\subseteq \bigcap_{i=1}^{n}
intZ(f_{i})\subseteq Z(f)$. But $\bigcap_{i=1}^{n}A_{i}$ contains
an element of $\cal{P}$, so $f\in C^{\cal{P}}(X)$. The proof of
the second equality is obvious.

(3) First, we observe that $C_{\infty\cal{P}}(X)$ is a proper
subset of $C(X)$. For, if $C_{\infty\cal{P}}(X)=C(X)$, then
$\emptyset\in\cal{P}$, which is a contradiction. On the other
hand, we have
\begin{center}
$\{x: |f(x)-g(x)|< \frac{1}{n}\}\supseteq \{x: |f(x)|<
\frac{1}{2n}\}\cap \{x: |g(x)|<
 \frac{1}{2n}\}$ and
\end{center}
\begin{center}
$\{x:|f(x)g(x)|<\frac{1}{n}\}\supseteq\{x: |f(x)|<
\frac{1}{\sqrt{n}}\}\cap \{x: |g(x)|<\frac{1}{\sqrt{n}}\}$.
\end{center}
Therefore, $C_{\infty\cal{P}}(X)$ is a proper subring of $C(X)$.
$\hfill\square$
\\\indent Recall that a family $\mathcal{P}$ of closed subsets of
$X$ is called an ideal of closed sets in $X$, if it satisfies in
the following conditions.
\begin{enumerate}
\item If $A, B\in\mathcal{P}$, then $A\cup B\in\mathcal{P}$.

\item If $A\in\mathcal{P}$ and $B\subseteq A$ with $B$ closed in
$X$, then $B\in\mathcal{P}$.
\end{enumerate}
In \cite{AK}, Acharyya and  Ghosh for ideal $\cal{P}$ of closed
subsets of $X$ defined $C_{\cal{P}}(X)$ and
$C_{\infty}^{\cal{P}}(X)$ as follows;
\begin{center}
 $C_{\cal{P}}(X)=\{f\in C(X):
cl(X\setminus Z(f))\in\cal{P}\}$; and
\end{center}
\begin{center}
$C_{\infty}^{\cal{P}}(X)=\{f\in C(X): \{x:
|f(x)|\geq\frac{1}{n}\}\in\cal{P}$, for each, $n\in\Bbb N\}.$
\end{center}
 In the next result we give a new presentation of these subrings.
 We note that for ideal $\mathcal{P}$ of closed sets, $C_{\infty}^{\cal{P}}(X)$ may be $C(X)$ but by Proposition \ref{2.2},
 $C_{\infty\cal{P}}(X)$ for each open filter $\mathcal{P}$ is a proper
 subring.
\begin{pro}
The following statements hold.
\begin{enumerate}
\item For every open filter base $\mathcal{P}$, there exists an
ideal $\mathcal{P}^{\prime}$ of closed sets such that
$C^{\cal{P}}(X)=C_{\cal{P}^{\prime}}(X)$ and
$C_{\infty\cal{P}}(X)=C_{\infty}^{\cal{P}^{\prime}}(X)$.

\item If $C_{\infty}^{\cal{P}}(X)$ is a proper subring of $C(X)$,
then there is an open filter base $\mathcal{Q}$ such that
$C_{\cal{P}}(X)=C^{\cal{Q}}(X)$ and
$C_{\infty}^{\cal{P}}(X)=C_{\infty\cal{Q}}(X)$.
\end{enumerate}
\end{pro}
{\noindent\bf Proof.} (1) Let $\mathcal{P}$ be an open filter
base. Consider $\mathcal{P}^{\prime}$ as follows;
\begin{center}
$\mathcal{P}^{\prime}=\{A:$ $A$ is closed and $A\subseteq
X\setminus B$ for some $B\in\cal{P}\}$.
\end{center}
 Then, it is easy to see that $\mathcal{P}^{\prime}$ is an ideal
of closed sets in $X$, $C^{\cal{P}}(X)=C_{\cal{P}^{\prime}}(X)$
and
$C_{\infty\cal{P}}(X)=C_{\infty}^{\cal{P}^{\prime}}(X)$.\\\indent
(2) Assume that $C_{\infty}^{\cal{P}}(X)$ is a proper subring of
$C(X)$ and $\cal{Q}$=$\{A\sub $X: X$\setminus A\in\cal{P}\}$. Then
we can see that $\cal{Q}$ is an open filter,
$C_{\cal{P}}(X)=C^{\cal{Q}}(X)$ and
$C_{\infty}^{\cal{P}}(X)=C_{\infty\cal{Q}}(X)$.$\hfill\square$
\begin{exa}
Let $X$ be a non-compact Hausdorff space and
$\cal{P}$=$\{A\subsetneq X: X\set A$ is compact$\}$. Then
$\cal{P}$ is an open filter base , $ C^{\cal{P}}(X)=C_{K}(X)$ and
$C_{\infty\cal{P}}(X)=C_{\infty}(X)$.
\end{exa}
\begin{exa}\label{a}
Let $\cal{P}$=\{$A$: $A$ is open and $X\setminus A$ is
Lindel\"{o}f$\}$  and $X$ be a non-Lindel\"{o}f space. Then
$\cal{P}$ is an open filter and we have
\begin{center}
$C_{\infty\cal{P}}(X)$=$\{f: X\set Z(f)$ is a Lindel\"{o}f subset
of $X\}$ and
\end{center}
\begin{center}
$C^{\cal{P}}(X)$=$\{f: \overline{X\set Z(f)}$ is a Lindel\"{o}f
subset of $X\}$.
\end{center}
 To see this, let $f\in C_{\infty\cal{P}}(X)$.
Then for each $n\in\Bbb N$, $\{x: |f(x)|\geq \frac{1}{n}\}\sub
X\set A$ for some $A\in\cal{P}$, so it is a Lindel\"{o}f subset of
$X$.
 On the other hand we have $X\set Z(f)=\bigcup_{n=1}^{\infty}\{x: |f(x)|\geq\frac{1}{n}\}$, hence $X\set Z(f)$ is a Lindel\"{o}f subset of $X$,
 by Theorem 3.8.5 in \cite{E}. If $X\set Z(f)$ be a Lindel\"{o}f subset of $X$, then $\{x: |f(x)|\geq\frac{1}{n}\}$ is a Lindel\"{o}f
subset of $X$ so $\{x: |f(x)|<\frac{1}{n}\}$ contains an element
of $\cal{P}$, i.e, $f\in C_{\infty\cal{P}}(X)$. Similarly we may
prove that $ C^{\cal{P}}(X)$=$\{f: \overline{X\set Z(f)}$ is a
Lindel\"{o}f subset of $X$\}.
\end{exa}
In the sequel, we assume that $\cal{P}$ is an open filter base for
a filter $\cal{F}$.
\begin{pro}
If the complement of any element of $\cal{P}$ is Lindel\"{o}f,
then $C_{\infty\cal{P}}(X)\sub \bigcap_{p\in \nu X\set X}M^{p}$.
Where by $\nu X$ we mean the real compactification of $X$ (see
\cite{G}).
\end{pro}
{\noindent\bf Proof.} Let $f\in C_{\infty\cal{P}}(X)$ and $p\in\nu
X\set X$ (i.e., $M^{p}$ is a free real maximal ideal). Then for
any $x\in X\set Z(f)$ there exist $f_{x}\in M^{P}$ such that $x\in
X\set Z(f_{x})$. Hence $X\set Z(f)\sub\bigcup_{x\in X} X\set
Z(f_{x})$. Now, by Example \ref{a} and hypothesis, $X\set Z(f)$ is
Lindelof, so there is a countable subset $S$ of $X$ such that
$X\setminus Z(f)\sub\bigcup_{x\in S}X\setminus Z(f_{x})$ and each
$f_{x}\in M^{p}$. On the other hand, there exists $h\in C(X)$ such
that $Z(h)=\bigcap_{x\in S}Z(f_{x})$. Therefore $h\in M^{p}$  and
$Z(h)\sub Z(f)$. But $M^{P}$ is a $z$-ideal, so $f\in M^{p}$,
i.e., $C_{\infty\cal{P}}(X)\sub \bigcap_{p\in \nu X\set
X}M^{p}$.$\hfill\square$
\\ \indent
Recall that, a subset $A$ of $X$ is called bounded (relative
pseudocompact) subset, if for every function $f\in C(X)$, $f(A)$
is a bounded subset of $\Bbb R$, see \cite{M}.
\begin{exa}
$\cal{P}$=\{$A$: $A$ is open and $X\setminus A$ is pseudocompact
$\}$ and $X$ be a non-pseudocompact space. Then $\cal{P}$ is an
open filter base and
\begin{center}
$C^{\cal{P}}(X)=C_{\psi}(X)=\{f:\overline{X\set Z(f)}$ is
pseudocompact$\}$, and
\end{center}
\begin{center}
$C_{\infty\cal{P}}(X)$=$\{f: \{x: |f(x)|\geq\frac{1}{n}\}$ is  pseudocompact$ \}$.
\end{center}
 For, suppose that $f\in  C^{\cal{P}}(X)$, then $Z(f)\supseteq A$ for some $A\in \cal{P}$. Hence $\overline{X\set Z(f)}\sub X\set A$. This
  implies that $\overline{X\set Z(f)}$ is a bounded subset. Now  by \cite[Theorem 2.1]{M}, $\overline{X\set Z(f)}$ is pseudocompact, i.e, $f\in C_{\psi}(X)$.
  If $f\in C_{\psi}(X)$, then $X\set \overline{X\set Z(f)}$ is an element of $\cal{P}$ and $Z(f)\supseteq X\set \overline{X\set Z(f)}$, i.e, $f\in  C^{\cal{P}}(X)$.
  For more details about $C_{\psi}(X)$, the reader referred to  \cite{M}. Similarly we may
prove the second equality.
 \end{exa}
\begin{Rem}
Let $\cal{P}$=$\{A: X\setminus A$ is finite$\}$ and $X$ is an
infinite space. Then $ C^{\cal{P}}(X)=C_{F}(X)$, and
$C_{\infty\cal{P}}(X)=\{f: \{x: |f(x)|\geq\frac{1}{n}\}$ is finite
for each $n\in\Bbb N\}$. In this case,
$C_{\infty\cal{P}}(X)=C_{F}(X)$ \ifif the set of isolated points
of $X$ is  finite. To see this, let $\{x_{1},x_{2},...x_{n},...\}$
be a  subset of isolated points in $X$. Define $f_{n}(x)=\left\{
\begin {array}{cc}
\frac{1}{n}& x=x_{n}\\
0& x\neq x_{n} \\
\end {array}
\right.$ and $f(x)=\sum_{n=0}^{\infty}f_{n}(x)$. Then $f\in C(X)$,
$f(x_{n})=\frac{1}{n}$ and $X\set Z(f)=\{x_{1},x_{2},...\}$, i.e.,
$f\notin C_{F}(X)$. On the other hand $\{x: |f(x)|<\frac{1}{n}\}$
contains  $X\set\{x_{1},x_{2},...,x_{n}\}$, so $f\in
C_{\infty\cal{P}}(X)$. This is a contradiction. Conversely, let
$f\in C_{\infty\cal{P}}(X)$ and $f\notin C_{F}(X)$. Then there is
$\{x_{1},x_{2},...x_{n},...\}$ such that $f(x_{n})\neq 0$, so
$|f(x_{n})|>\frac{1}{k_{n}}$ for some $k_{n}\in\N$, i.e.,
$x_{n}\in \{x: |f(x)|\geq \frac{1}{k_{n}}\}$, but $ \{x: |f(x)|>
\frac{1}{k_{n}}\}$ is a finite open set so each $x_{n}$ is an
isolated point, this is a contradiction.
\end{Rem}
\begin{defn}
A topological space $X$ is called a $\textit{local}$
$\textit{space}$ provided that, there exists an open filter base
$\cal{P}$ for a filter $\cal{F}$ on space $X$ where $\cal{F}$ has
no cluster point (i.e.,
$\bigcap_{A\in\cal{P}}\overline{A}=\varnothing$).
\end{defn}
\begin{exa}
Any locally compact non-compact space $X$ is a local space. For
see this, let $\cal{P}$=$\{A\sub X: X\set A$ is compact$\}$. Then
if $x\in \bigcap_{A\in \cal{P}}\overline{A}$, by locally
compactness of $X$, there exists a compact subset $K\subsetneqq X$
such that $x\in int(K)$, but $X\set K$ is in $\cal{P}$, hence
$x\in \overline{X\set K}$, i.e., $int(K)\cap (X\set
K)\neq\emptyset$, this is a contradiction, thus $\bigcap_{A\in
\cal{P}}\overline{A}=\tohi$.
\end{exa}
\begin{exa}
\cite[4. M]{G} Let $X$ be an uncountable space in which all points
are isolated points except for a distinguished point $s$. A
neighborhood of $s$ is any set
 containing $s$ which complement is countable. Then $X$ is a local space. To see this,
 let $Y=\{x_{1},x_{2},...\}$ be a countable subset of $X$, where $s\notin Y$. Put $A_{n}=\{x_{n},x_{n+1},...\}$ and $\cal{P}$=$\{A_{n}: n\in\Bbb N\}$.
 Then $\cal{P}$ is an open  filter base on $X$. Now, for any $n\in \Bbb N$, the set $X\set A_{n}$  is a neighborhood of
 $s$,
so $s\notin \bigcap_{A_{\widetilde{}n}\in
\cal{P}}\overline{A_{n}}$, thus $\bigcap_{A\in
\cal{P}}\overline{A}=\tohi$.
\end{exa}
 In the following we see an example of a topological space which is not a local space.
\begin{exa}\cite[4. N]{G}
For each $n\in\N$, let $A_{n}=\{n,n+1,...\}$ and $E=\{A_{n}$:
$n\in\N\}$. Then $E$ is a base for a free ultrafilter say
$\cal{E}$ on $\N$. Let $X=\N\cup \{\sigma\}$ which points in $\N$
are isolated point and a neighborhood of $\sigma$ is of the form
$U\cup\{\sigma\}$ which $U\in\cal{E}$. Note that any set contains
$\sigma$ is closed. Now if there is an open base  $\cal{P}$ for
some filter $\cal{F}$ on $X$  such that $\cal{F}$ has no cluster
point, then  there exists $F\in\cal{F}$ such that
$\sigma\notin\overline{F}$, but $\sigma$ has a neighborhood say
$U\cup\{\sigma\}$ such that  $U\in\cal{E}$ and
$U\cup\{\sigma\}\sub X\set\overline{F}$.
 Since $E$ is a base for $\cal{E}$, then there exists $n\in\N$ and $A_{n}\in E$ such that $A_{n}\cup\{\sigma\}=\{n,n+1,...\}\cup\{\sigma\}\sub U\cup\{\sigma\}\sub X\set\overline{F}$.
 On the other hand, for points $x=1,2,...n-1$ there exist $F_{1},F_{2},...F_{n-1}$ in $\cal{F}$, such that $i\in X\set F_{i}$ for
  $1\leq i\leq n$, so $X\sub (X\set\overline{F})\cup\bigcup_{i=1}^{n-1}(X\set F_{i})$, therefore $\emptyset\in\cal{F}$. This a contradiction, i.e., $X$ is not a
  local space.
\end{exa}
We have already observe that $ C^{P}(X)$ is a $z$-ideal in $C(X)$.
In the following proposition we find a condition over which
$C^{\cal{P}}(X)$ is a free ideal.
\begin{pro}\label{2.3}
Let $\cal{P}$ be an open  filter base for a  filter $\cal{F}$.
Then
 $ C^{\cal{P}}(X)$ is a free ideal \ifif $\cal{F}$ has no cluster point (i.e., $\bigcap_{A\in\cal{P}}\overline{A}=\varnothing$).
\end{pro}
{\noindent\bf Proof.} Let $\bigcap_{F\in\cal{F}}\overline{F}\neq
\emptyset.$ Then there exists $x\in
\bigcap_{F\in\cal{F}}\overline{F}.$ By hypothesis, there is $f\in
C^{\cal{P}}(X)$ such that $x\in X\setminus Z(f)$. On the other
hand there is $A\in\cal{F}$ such that $X\setminus Z(f)\subseteq
X\setminus A$. Hence $(X\setminus Z(f))\cap A=\emptyset$. But
$x\in \overline{A}$ implies that $(X\setminus Z(f))\cap
A\neq\emptyset$. This is a contradiction. Conversely, let $x\in
X$. We have $\bigcap_{F\in\cal{F}}\overline{F}= \emptyset$, so
there
 exists $F\in\cal{F}$ such that $x\notin \overline{F}$. By completely regularity of $X$, there exists $f\in C(X)$ such that $f(x)=1$ and $f(\overline{F})=0$.
 Hence $f\in C^{\cal{P}}(X)$ and $x\notin Z(f)$, i.e., $ C^{P}(X)$ is a free
 ideal. $\hfill\square$
\\ \indent
It is easy to see that $X$ is a locally compact non-compact space
\ifif $\cal{P}$=$\{A\subsetneq X: X\set A$ is compact$\}$ is an
open filter base for some filter $\mathcal{F}$
 with no cluster point. So, by the above proposition we have the following corollaries.
\begin{cor}
$C_{K}(X)$ is a free ideal \ifif  $X$ is a locally compact non-compact space.
\end{cor}
\begin{cor}
 A space $X$ is a local space \ifif $ C^{\cal{P}}(X)$ is a free ideal for some open filter base $\cal{P}$ on $X$.
 \end{cor}
 {\noindent\bf Proof.} By Proposition \ref{2.3}, the verification is
 immediate. $\hfill\square$
\begin{pro}\label{0}
Let $\cal{P}$ be an open filter base. The following statements are
equivalent.
\begin{enumerate}
\item Every element of $\cal{P}$ is dense in $X$.

\item $C_{\infty\cal{P}}(X)=(0).$

\item $ C^{\cal{P}}(X)=(0).$
\end{enumerate}
\end{pro}
{\noindent\bf Proof.} (1)$\Rightarrow$(2) Let for every
$A\in\cal{P}$, $\overline{A}=X$ and $f\in C_{\infty\cal{P}}(X)$.
Then the set $\{x: |f|\leq\frac{1}{n}\}=X$, so for any $n>1$ we
have $\{x: |f|<\frac{1}{n-1}\}\supseteq\{x:
|f|\leq\frac{1}{n}\}=X$, i.e., $f=0$.

(2)$\Rightarrow$(3) This is evident.

 (3)$\Rightarrow$(1) Suppose that  $ C^{\cal{P}}(X)=0$ and $A\in \cal{P}$. If $\overline{A}\neq X$, then there
 exists
 $x\in X\set \overline{A}$, hence we define $f\in C(X)$ such that $f(x)=1$ and $f(\overline{A})=0$, i.e., $f\in C_{\cal{P}}(X)=0$, which implies that
 $f=0$, this is a contradiction. $\hfill\square$
\begin{cor}
Let $X=\Bbb Q$ with usual topology and $\cal{P}$=$\{A\subset \Bbb
Q: \Bbb Q\set A$ is compact $\}$. Then
$C_{\infty\cal{P}}(X)=C_{\infty}(X)=(0)$.
\end{cor}
{\noindent\bf Proof.} Every element of $\cal{P}$ is dense in $X$,
so by Proposition \ref{0},
$C_{\infty\cal{P}}(X)=C_{\infty}(X)=(0)$. $\hfill\square$
\begin{defn}
Let $\cal{F}$ be a filter on $X$. A subset $A$ of $X$ is an
$\cal{F}$-$CG_{\delta}$ if  $A=\bigcap_{i=1}^{\infty}A_{i}$, where
each $A_{i}\in\cal{F}$ and is open  and  for each $i$, $X\set
A_{i}$ and $\overline{A_{i+1}}$ are completely separated (see
\cite{G}).
\end{defn}
\begin{exa}
Let $\cal{F}$=$\{F\subsetneqq X: X\setminus F$ is compact$\}$ and
$X$ is a non-compact Hausdorff  space. Then for every open locally
compact $\sigma $-compact subspace $A$, $X\setminus A$ is an
$\cal{F}$-$CG_{\delta}$. Since, by \cite[3. p. 250]{E},
$A=\bigcup_{i=1}^{\infty}A_{i}$ where $A_{i}\subseteq intA_{i+1}$
and each $A_{i}$ is compact so $X\setminus A$ is an
$\cal{F}$-$CG_{\delta}$ set.
\end{exa}
In the following lemma, we give a  characterization of a closed
$\cal{F}$-$CG_{\delta}$ subset.
\begin{lem}\label{2.4} Let $A$ be a closed
subset of space  $X$. Then $A$ is an $\cal{F}$-$CG_{\delta}$ set
 \ifif $A=Z(f)$ for some $f\in
C_{\infty\cal{P}}(X)$.
\end{lem}
{\noindent\bf Proof.}  Let $A$ be an $\cal{F}$-$CG_{\delta}$. Then
$A=\bigcap_{n=1}^{\infty}A_{i}$, where each $A_{n}$ is a an
element of $\cal{F}$, $X\set A_{n}$ and $\overline{A_{n+1}}$ are
completely separated. Now  for each $n\in \Bbb N$, there exists
$f_{n}\in C(X)$ such that $f_{n}(\overline{A_{n+1}})=0$,
$f_{n}(X\set A_{n})=1$, then $f=\sum\frac{1}{{2^n}}f_{n}$ is an
element of $C(X)$, by Weierstrass M-test. Clearly $A= Z(f)$. We
claim that $f\in C_{\infty\cal{P}}(X)$. Let $x_{0}\in A_{n+1}$.
Then $f_{1}(x_{0})=f_{2}(x_{0})=...f_{n}(x_{0})=0$ and so
$f(x_{0})\leq \frac{1}{2^{n+1}}+\frac{1}{2^{n+1}}+...\leq
\frac{1}{2^{n}}< \frac{1}{n}$. Therefore $x_{0}\in \{x: |f(x)<
\frac{1}{n}\}$, and hence $A_{n+1}\sub\{x: |f(x)<
 \frac{1}{n}\} $, i.e, $f\in
 C_{\infty\cal{P}}(X)$. Conversely, suppose that $A=Z(f)$ for some $f\in C_{\infty\cal{P}}(X)$. Then $A=\bigcap_{n=1}^{\infty}A_{n}$, where
 $A_{n}=\{x: |f(x)|< \frac{1}{n}\}\in\cal{F}$ for each $n\in\Bbb
N$, $X\setminus A_{n}$ and $\overline{A_{n+1}}$ are disjoint
zero-sets, and hence completely separated, i.e., $A$ is an
$\cal{F}$-$CG_{\delta}$. $\hfill\square$

 \begin{pro}\label{2.5}
 $C_{\infty\cal{P}}(X)=C^{\cal{P}}(X)$ if and only if every closed $\cal{F}$-$CG_{\delta}$  is an element of $\cal{F}$.
 \end{pro}
{\noindent\bf Proof.}
 Suppose that condition holds.  We know that  $C^{P}(X)$ is a subset of
 $C_{\infty\cal{P}}(X)$. It is enough to prove that
$C_{\infty\cal{P}}(X)\subseteq  C^{\cal{P}}(X)$. Let $f\in
C_{\infty\cal{P}}(X)$, then by Lemma \ref{2.3}, $ Z(f)$ is a
closed $\cal{F}$-$CG_{\delta}$. Hence $ Z(f)$ contains
 an element of $\cal{P}$ , i.e, $f\in  C^{\cal{P}}(X)$. Conversely, suppose that
$C_{\infty\cal{P}}(X)=C_{\cal{P}}(X)$ and $A$ is a closed
$\cal{F}$-$CG_{\delta}$. By lemma \ref{2.3}, $A= Z(f)$ for some
$f\in
 C_{\infty\cal{P}}(X)$, now $f\in
  C^{\cal{P}}(X)$ implies that $A=Z(f)$
  contains an element of $\cal{P}$, i.e, $A\in\cal{F}$. $\hfill\square$
\\ \indent
In  the above proposition we have seen that if every closed
$\cal{F}$-$CG_{\delta}$  is an element of $\cal{F}$, then
$C_{\infty\cal{P}}(X)$ is an ideal of $C(X)$. But in general,
  $C_{\infty\cal{P}}(X)$ may not be an ideal of $C(X)$ as we will see in the sequel.
\begin{exa}
Let $\cal{P}$=$\{\RE\set[\frac{1}{n}, n]$:  $n\in \Bbb N \}$. Then
it is easy to see that $\cal{P}$ is an open filter base on $\RE$.
Now, we show that $C_{\infty\cal{P}}(\RE)\subsetneqq
C_{\infty}(\RE)$ is not an ideal of $C(\RE)$. For see this,
Consider
\begin{center}
 $f(x)=\left\{
\begin {array}{cc}
0& x\leq 0 \\
x& 0\leq x\leq 1 \\
\frac{1}{x^{2}}&1\leq x \\
\end {array}
\right.$,
$g(x)=\left\{
\begin {array}{cc}
0& x\leq 0 \\
x& 0\leq x\leq 1 \\
x^{2}&1\leq x\\
\end {array}
\right..$
\end{center}
 Then $f\in C_{\infty\cal{P}}(\RE)$, $g\in C(\RE)$ and we
have
\begin{center}
$(fg)(x)=\left\{
\begin {array}{cc}
0& x\leq 0 \\
x^{2}& 0\leq x\leq 1 \\
1&1\leq x\\
\end {array}
\right.$
\end{center}
But $fg\notin C_{\infty\cal{P}}(\RE)$, because $\{x: |(fg)(x)|<
1\}$=$(-\infty,1)$. On the other hand $\frac{1}{x^{2}+1}\in
C_{\infty}(\RE)$ and is not in $C_{\infty\cal{P}}(\RE)$, so
$C_{\infty\cal{P}}(\RE)\subsetneqq C_{\infty}(\RE)$.
\end{exa}
In the following, we provide  an example of open filter base
$\cal{P}$ for which $C_{\infty\cal{P}}(X)$ is an ideal of $C(X)$.
\begin{exa}
$\cal{P}$=\{$A$: $A$ is open and $X\setminus A$ is
Lindel\"{o}f$\}$  and $X$ be a non-Lindel\"{o}f space. Then by
Example \ref{a}, $C_{\infty\cal{P}}(X)$=$\{f:X\set Z(f)$ is
Lindel\"{o}f subset of $X\}$, which is an ideal of $C(X)$. For,
let $f,g\in
 C_{\infty\cal{P}}(X)$, then $X\set Z(f+g)\subseteq (X\set Z(f))\cup
 (X\set Z(g))$ implies $\{x: |(f+g)(x)|<\frac{1}{n}\}\supseteq X\set((X\set Z(f))\cup
 (X\set Z(g))$, which contains  an element of $\cal{P}$, i.e, $f+g\in C_{\infty\cal{P}}(X)$. If $f\in C_{\infty\cal{P}}(X)$ and $g\in C(X)$, then $X\set
Z(fg)\subseteq X\set Z(f)$ implies that $\{x: |(f.g)(x)|<\frac{1}{n}\}$  contains  an element of $\cal{P}$, i.e, $f.g\in C_{\infty\cal{P}}(X)$.
\end{exa}
Azarpanah and  Soundarajan in \cite{AS} have found some equivalent
conditions for which $C_{\infty}(X)$ is an ideal of $C(X)$ (i.e.,
$\cal{P}$=$\{A\subsetneqq X: X\setminus A$ is compact$\}$). In the
following theorem, we find some equivalent conditions for a larger
class of $\cal{P}$  which $C_{\infty\cal{P}}(X)$ is an ideal of
$C(X)$ and by this theorem give several corollary.
\begin{thm}\label{2.6}
The following statements are
equivalent.
\begin{enumerate}
\item $C_{\infty\cal{P}}(X)$ is an ideal of $C(X)$ and for any
$F\in \cal{F}$, $X\set \overline{F}$ is bounded.

\item  The set of non-cluster points of filter $\cal{F}$ is
bounded.

\item  The complement of every closed $\cal{F}$-$CG_{\delta}$ is
bounded.
\end{enumerate}
\end{thm}

{\noindent\bf Proof.} (1)$\Rightarrow$(2). Let $A$ be the set of
non-cluster points of filter $\cal{F}$. If $A$ is not bounded,
 then there exist $h\in C(X)$ and a discrete subset $C=\{x_{1},x_{2},x_{3},...\}\sub A$ such that for each $n\in \Bbb
N$, $h(x_{n})\geq n$.  $A$ is open so any $\{x_{n}\}$ is an isolated point, therefore we can define
$f_{n}(x)=\left\{
\begin {array}{cc}
\frac{1}{n}& x=n \\
0& x\neq n \\
\end {array}
\right.$ and $f(x)=\sum_{n=0}^{n=\infty}f_{n}(x)$ such that
$f_{n}\in C(X)$ and so $f\in C(X)$. We have
$\{x:|f|<\frac{1}{n}\}=\{x_{n+1},x_{n+2},...\}$. Now any $x_{n}\in
X\set \overline{F_{n}}$ for some $F_{n}\in \cal{F}$. This implies
that $X\set \{x_{1},...,x_{n}\}$ contains an element of $\cal{F}$
hence contains an element of $\cal{P}$, i.e., $f\in
C_{\infty\cal{P}}(X)$. But we have $\{x:|fh|<\frac{1}{n}\}=X\set
\{x_{1},x_{2}...,\}$, if there is $P\in\cal{P}$ such that $X\set
\{x_{1},x_{2}...,\}\supseteq P$, then $\{x_{1},x_{2},...\}\sub
X\set \overline{P}$, which contradict by hypothesis, so $fh\notin
C_{\infty\cal{P}}(X)$, i.e., $C_{\infty\cal{P}}(X)$ is not an
ideal of $C(X)$, this is a contradiction.

(2)$\Rightarrow$(3). It is easily seen that the  complement of
every closed $\cal{F}$-$CG_{\delta}$  is a subset of  non-cluster
points of the filter $\cal{F}$ so is bounded.

(3)$\Rightarrow$(1). By Lemma \ref{2.2}, $C_{\infty\cal{P}}(X)$ is
a subring of $C(X)$, it is enough to prove that $fg\in
C_{\infty\cal{P}}(X)$ for every $f\in C(X)$ and $g\in
C_{\infty\cal{P}}(X)$. By Lemma \ref{2.3}, $ Z(g)$ is an
$\cal{F}$-$CG_{\delta}$ and hence, by (3), $Y=X\set Z(g)$ is a
bounded subset of $X$, so $f(Y)$ is a bounded subset of $\Bbb R$.
Now it is easy to see that $g^{\frac {1}{3}}\in
C_{\infty\cal{P}}(X)$, since $g\in C_{\infty\cal{P}}(X)$, moreover
$(f(g)^{\frac{1}{3}})(X)= (f(g)^{\frac{1}{3}})(Y)\cup \{0\}$.
Since $f(Y)$ is a bounded subset of $\Bbb R$ and $g^{\frac
{1}{3}}\in C_{\infty\cal{P}}(X)$ is a bounded function on $X$
(note that $X\set Z(g)$ is bounded), we get
$(f(g)^{\frac{1}{3}})(Y)$ is a bounded set in $\Bbb R$, this
implies that $(f(g)^{\frac{1}{3}})$ is a bonded on $X$ and so
belong to $C^{*}(X)$. Since $C_{\infty\cal{P}}(X)$ is a ring,
$g^{\frac{2}{3}}\in C_{\infty\cal{P}}(X)$. However, since for any
$f\in C_{\infty\cal{P}}(X)$, $X\setminus Z(f)$ is bounded, then
$C_{\infty\cal{P}}(X)$ is an ideal of $C^{*}(X)$. Therefore $fg=
(f(g)^{\frac{1}{3}})(g^{\frac{2}{3}})\in C_{\infty\cal{P}}(X)$,
thus $C_{\infty\cal{P}}(X)$ is an ideal of $C(X)$. Now let
$F\in\cal{F}$ and $X\set\overline{F}$ is unbounded. Then it
contains infinity set of  isolated points say
$A=\{x_{1},x_{2},...\}\sub (X\set \overline{F})$, such that $A$ is
unbounded. We have $X\set A =\bigcap_{i=1}^{\infty}B_{i}$ where
$B_{i}=X\set\{x_{1},...,x_{i}\}$. Clearly for each $i$
$B_{i}\in\cal{F}$ and $X\set B_{i}$, $\overline{B_{i+1}}$  are
completely separated. Hence $X\set A$ is an
$\cal{F}$-$CG_{\delta}$, by hypothesis, $X\set(X\set A)=A$ is
bounded, this is a contradiction. This completes the proof.
$\hfill\square$
\begin{cor}\label{2.7}
If $X$ is a  pseudocompact  space, then for
any open filter base $\cal{P}$, $C_{\infty\cal{P}}(X)$ is an
ideal of $C(X)$.
\end{cor}
{\noindent\bf Proof.} If $X$ is a completely regular pseudocompact
Hausdorff  space and $\cal{P}$ be an open filter base for filter
$\cal{F}$, then any subset of $X$ is bounded so the set of
non-cluster point of $\cal{F}$ is bounded, thus by Theorem
\ref{2.6}, $C_{\infty\cal{P}}(X)$ is an ideal of $C(X)$.
$\hfill\square$
\begin{cor}\label{2.9}
Let $X$ be a local space. Then for any open filter base $\cal{P}$, $C_{\infty\cal{P}}(X)$ is an ideal of $C(X)$  and for any $A\in \cal{P}$,
$X\set \overline{A}$ is bounded \ifif $X$ is a pseudocompact non-compact space.
\end{cor}
{\noindent\bf Proof.} If $X$ is pseudocompact, then by Corollary
\ref{2.7}, for any open filter base $\cal{P}$,
$C_{\infty\cal{P}}(X)$ is an ideal of $C(X)$ and for any $A\in
\cal{P}$, $X\set \overline{A}$ is bounded. Now let $X$ be a local
space, then  there exist an open filter base $\cal{P}$ for some
filter $\cal{F}$ on $X$ such that $\cal{F}$ has no cluster point
so $X$ is the set of non-cluster point of filter $\cal{F}$. Hence
by Theorem \ref{2.6}, $X$ is bounded, i.e., $X$ is pseudocompact.
$\hfill\square$
\begin{cor}
 Let $X$ be a non-pseudocompact space and $\cal{P}$=$\{A$: $A$ is open and $X\setminus A$ is bounded $\}$. Then $C_{\infty\cal{P}}(X)$ is an ideal of
 $C(X)$ \ifif any union of the interior of closed bounded subsets is a bounded subset.
\end{cor}
{\noindent\bf Proof.} If any union of the interior of closed
bounded subsets is a bounded subset, then the set of non-cluster
points of open filter $\cal{P}$ is bounded so by Theorem
\ref{2.6}, $C_{\infty\cal{P}}(X)$ is an ideal of $C(X)$.
Conversely, If $A=\bigcup_{\alpha\in S}$int$A_{\alpha}$ where for
each $\alpha\in S$, $A_{\alpha}$ is a closed bounded set, then we
have $X\setminus A_{\alpha}\in\cal{P}$ and
int$A_{\alpha}=X\setminus \overline{(X\setminus A_{\alpha})}$, so
$A$ is contained in the set of non-cluster points of open filter
$\cal{P}$, i.e., $A$ is bounded. $\hfill\square$
\begin{cor}\label{i}
Let $X$ be an an infinite space and  $\cal{P}$=$\{A\subsetneq X$: $X\setminus A$ is finite $\}$. Then $C_{\infty\cal{P}}(X)$ is an ideal of $C(X)$ \ifif the set of isolated points of $X$ is bounded.
\end{cor}
{\noindent\bf Proof.} Let $A\in\cal{P}$. We have $X\setminus
\overline{A}$ is an open finite subset, thus the set of
non-cluster points of $\cal{P}$ is contained in the set  of
isolated points of $X$, so is bounded, hence by Theorem \ref{2.6},
$C_{\infty\cal{P}}(X)$ is an ideal of $C(X)$. Conversely, let $A$
be the set of isolated points of $X$. Then $A=\bigcup_{x\in
A}\{x\}$,
 each $\{x\}$ is a clopen subset, $X\setminus \{x\}\in\cal{P}$ and $\{x\}$=int$\{x\}=X\setminus \overline{(X\setminus \{x\})}$, so $A$ is contained in the set of
 non-cluster points of open filter $\cal{P}$, i.e, $A$ is bounded.
 $\hfill\square$
\begin{exa}
(a). If $\cal{P}$=$\{A\subsetneq \Bbb R$: $\Bbb R\setminus A$ is bounded $\}$. Then $C_{\infty\cal{P}}(\Bbb R)$ is not an ideal of $C(\Bbb R)$. Because $\bigcup_{n=1}^{\infty}(0, n)$ is not bounded.
\\ \indent (b). If $\cal{P}$=$\{A\subsetneq \Bbb R$: $\Bbb R\setminus A$ is finite $\}$. Then by Corollary \ref{i}, $C_{\infty\cal{P}}(\Bbb R)$ is an ideal of $C(\Bbb R)$.
\end{exa}
\begin{Rem}\label{2.10} Any closed bounded in a normal space
 is a pseudocompact  and any pseudocompact  Lindel\"{o}f space is
 compact, so  if $X$ is a realcompact normal
space and $\cal{P}$ equal be the set of all subsets whose complements are bounded subsets of
$X$, then $C_{\infty\cal{P}}(X)=C_{\infty}(X)$, for example  in $\Bbb R$  if $\cal{P}$ equal be the set of all  subsets of $\Bbb R$, whose complements are bounded,
then $C_{\infty\cal{P}}(\Bbb R)=C_{\infty}(\Bbb R)$. In particularly,
if $X$ is a Lindel\"{o}f space and $\cal{P}$ equal be the set of all  subsets of $X$, whose complements are bounded,
then $C_{\infty\cal{P}}(X)=C_{\infty}(X)$.
\end{Rem}



\section {$C^{\cal{P}}(X)$ as an essential ideal.}
Topological spaces $X$ for which $C_{\infty}(X)$ (resp.,
$C_{K}(X)$) is an essential ideal was characterized by Azarpanah,
in \cite{AZ}. In this section we characterize topological spaces
$X$ for which $C^{\cal{P}}(X)$ is an essential ideal.
\begin{pro}\label{2.10}
An ideal $E$ in $C_{\infty\cal{P}}(X)$ is an essential ideal \ifif
$\bigcap Z[E]$ does not contain a subset $V$ where $X\setminus
V\in\mathcal{F}$ and $intV\neq\emptyset$.
\end{pro}
{\noindent\bf Proof.} Let $X\setminus V=F\in\cal{F}$,
$\bigcup_{f\in E} coz(f)\subseteq F$ and $intV\neq\emptyset$,
i.e., $\overline{F}\neq X$. Then there exist $x\in X$ such that
$x\notin \overline{F}$, it follows that there is $f\in C(X)$ such
that $f(x)=1$, $f(F)=0$, i.e., $f\in C_{\infty\cal{P}}(X)$. Now
for any $g\in E$ we have $X\set Z(f)\sub X\set F\sub  Z(g)$, i.e.,
$fg=0$ so $(f)\cap E=0$, which contradicts the essentiality of
$E$. Conversely, let $0\neq f\in C_{\infty\cal{P}}(X)$. Then there
is $a\in X$ such that $|f(a)|>\frac{1}{n}$ for some $n\in \Bbb N$,
hence $a\in X\set \{x: |f|\leq\frac{1}{n}\}$, i.e., $\{x:
|f|\leq\frac{1}{n}\}\neq X$. We know that $\{x:
|f|<\frac{1}{n}\}\in\cal{F}$. By hypothesis, $X\set \{x:
|f|<\frac{1}{n}\}\nsubseteq \bigcap Z[E]$. Therefore, there exists
$b\in X\set \{x: |f|\leq\frac{1}{n}\}$ and $g\in E$ such that
$g(b)\neq 0$, i.e, $fg\neq 0$ thus $E$ is an essential ideal in
$C_{\infty\cal{P}}(X)$. $\hfill\square$
\\\indent Recall that a collection $\cal{B}$ of open sets in a topological space
$X$ is called a $\pi$-base if every open set of $X$ contains a
member of $\cal{B}$. The reader is referred to \cite{BH},
\cite{BM}, \cite{E}, \cite{KM}, and \cite{MC}. The next result is
a generalization of  Theorem 3.2 in \cite{AZ}.
\begin{thm}\label{2.11}
$C^{\cal{P}}(X)$ is an essential ideal \ifif  $\{V:$ $V$ is open
and $X\setminus V\in\mathcal{F}\}$ is a $\pi$-base for $X$.
\end{thm}
{\noindent\bf Proof.}  Let $U$ be a proper open set in $X$. By
regularity of $X$, there exist a non-empty open set $V$ such hat
$V\sub clV\sub U$. Now find $f\in C(X)$ where $f(clV)=\{1\},
f(x)=0$ for some $x\notin U$. If $ X\set V\in\cal{F}$, there is
noting to proved. Suppose $ X\set V\notin\cal{F}$. If $V\sub Z(h)$
for every $h\in C^{\cal{P}}(X)$, then $V\sub \bigcap
Z[C^{\cal{P}}(X)]$, which implies that  $C^{\cal{P}}(X)$ is not an
essential ideal, by \cite[Theorem 3.1]{A}. Therefor there is some
$h\in C^{\cal{P}}(X)$ such that $V\cap (X\set Z(h))\neq \tohi$,
i.e, there is some $x_{0}\in V$ for which $h(x_{0})\neq 0$.
Clearly $fh\in C^{\cal{P}}(X)$. So $W=X\set Z(fh)$ is contained in
$X\set F$ for some $F\in\cal{F}$. If $W^{\prime}=W\cap V$, then
$W^{\prime}$ is a non-empty open set in $U$ and  $X\set
W^{\prime}\in\cal{F}$.

 Conversely,  We will prove that for every non-unit $g\in
C(X)$, $ C^{\cal{P}}(X)\cap (g)\neq 0$. Since $X\set Z(g)$ is an
open set, then there is an open set $U$ where $U\sub cl U\sub
X\set Z(g)$, and there is an open set $V\sub U$ such that $ X\set
V\in\cal{F}$. Then $V\sub U\sub X\set Z(g)$. Define $f\in C(X)$
such that $f(X\set V)=0 , f(x)=1$ for some $x\in V$. Since $ X\set
V\sub Z(f)$ so $f\in  C^{\cal{P}}(X)$ . On the other hand
$Z(g)\sub X\set V\sub Z(f)$ so $fg\neq 0$ and $fg\in
C_{\cal{P}}(X)\cap (g)$. $\hfill\square$

\section{ $C_{\infty\cal{P}}(X)$ as a $z$-ideal and a regular ring.}

We know that $C_{\infty\cal{P}}(X)$ is a subring of $C(X)$, in
this section, we see that $C_{\infty\cal{P}}(X)$ is a $z$-ideal
\ifif  every cozero-set  containing a $\cal{F}$-$CG_{\delta}$ is
an element of $\cal{F}$. Also we prove that,
$C_{\infty\cal{P}}(X)$ is a regular ring  (i.e., for each $f\in
C_{\infty\cal{P}}(X)$ there exists $g\in C_{\infty\cal{P}}(X)$
such that $f=f^{2}g$) \ifif every closed $\cal{F}$-$CG_{\delta}$
is an open subset and belong to $\cal{F}$.
\begin{pro}
The subring $C_{\infty\cal{P}}(X)$ is a $z$-ideal of $C(X)$ \ifif
every cozero-set  containing a closed  $\cal{F}$-$CG_{\delta}$ is
an element of $\cal{F}$.
\end{pro}
 {\noindent\bf Proof.} First, we prove that $Z(f)\subseteq Z(g)$ and $f\in C_{\infty\cal{P}}(X)$, implies that $g\in
 C_{\infty\cal{P}}(X)$. To see this, we know that
   $Z(f)\subseteq Z(g)\subseteq \{x: |g(x)|<\frac{1}{n}\}$, for all $n\in\Bbb N$. But $\{x: |g(x)|<\frac{1}{n}\}$ is a cozero-set and $Z(f)$ is a closed
   $\cal{F}$-$CG_{\delta}$. So, by hypothesis, $\{x: |g(x)|<\frac{1}{n}\}$ is an element of $\cal{F}$, i.e., $g\in C_{\infty\cal{P}}(X)$.
   Now, suppose that $f\in C_{\infty\cal{P}}(X)$ and $g\in C(X)$. Then
$Z(f)\sub Z(fg)$, shows that $fg\in C_{\infty\cal{P}}(X)$. Thus
$C_{\infty\cal{P}}(X)$ is a $z$-ideal of $C(X)$. Conversely,
Suppose that $X\setminus Z(f)$ is a cozero-set contains a closed
$\mathcal{F}CG_{\delta}$ subset $A$. By Lemma \ref{2.4},
 there exist $g\in C_{\infty\cal{P}}(X)$, such that $A=Z(g)$, so $Z(g)\subseteq X\setminus Z(f)$. Now we define $h=\frac{g^{2}}{f^{2}+g^{2}}$.
 Then $h\in C(X)$ and $Z(g)\subseteq Z(h)$, therefore $h\in C_{\infty\cal{P}}(X)$. On the other hand, for each $n\in\Bbb N$, $\{x: |h(x)|<\frac{1}{n}\}\subseteq X\setminus Z(f)$,
  hence $X\setminus Z(f)\in\cal{F}$. This completes the proof. $\hfill\square$
\begin{thm}\label{2.13}
$C_{\infty\cal{P}}(X)$ is a regular ring \ifif every closed $\cal{F}$-$CG_{\delta}$ is an open subset and belongs to $\cal{F}$.
\end{thm}
{\noindent\bf Proof.} First, we prove that every closed
$\cal{F}$-$CG_{\delta}$ is an open subset. By Lemma \ref{2.4},
every closed $\cal{F}$-$CG_{\delta}$ is of the form $Z(f)$ for
some $f\in C_{\infty\cal{P}}(X)$. But $Z(f)=Z(f\wedge n)$, for
each $n\in \Bbb N$ and $\{x: |f\wedge n|<\frac{1}{m}\}=\{x:
|f|<\frac{1}{m}\}$. So we can let $f$ is bounded. Regularity of
$C_{\infty\cal{P}}(X)$ implies that, there exists $g\in
C_{\infty\cal{P}}(X)$ such that $f=f^{2}g$. Then $X\setminus
Z(1-fg)\subseteq intZ(f)$. If $x\in Z(f)\setminus intZ(f)$, then
$x\in Z(1-fg)$, which contradict $x\in Z(f)$, i.e., $Z(f)$ is an
open subset. On the other hand for every $x\in X\setminus Z(f)$,
$g(x)=\frac{1}{f(x)}$ and hence $g(x)\geq \frac{1}{n}$, where $n$
is an upper bounded for $|f|$. Therefore $X\setminus Z(f)\subseteq
\{x: |g|\geq \frac{1}{n}\}$, i.e., $ Z(f)\supseteq \{x: |g|
<\frac{1}{n}\}$. But $\{x: |g| <\frac{1}{n}\}$ contains an element
of $\cal{P}$ so $Z(f)\in\cal{F}$. Conversely, Suppose $f\in
C_{\infty\cal{P}}(X)$. $Z(f)$ is a closed $\cal{F}$-$CG_{\delta}$
so by hypothesis, is an open subset which belong to $\cal{F}$. We
define $g(x)=0$ for $x\in Z(f)$ and $g(x)=\frac{1}{f(x)}$ for
$x\in X\setminus Z(f)$. Then $g\in C(X)$, $f=f^{2}g$ and $\{x: |g|
<\frac{1}{n}\}\supseteq Z(f)$. , i.e., $g\in
C_{\infty\cal{P}}(X)$. $\hfill\square$

\begin{cor}
(a) Let $\cal{P}$=\{$A$: $A$ is open and $X\setminus A$ is
Lindel\"{o}f$\}$ and $X$ is a non-Lindel\"{o}f space. Then
$\cal{P}$ is an open filter and $C_{\infty\cal{P}}(X)$ is a
regular ring \ifif every closed $\cal{P}$-$CG_{\delta}$ is an open
subset.

(b) $C_{\infty}(X)$ is a regular ring \ifif every open locally
compact $\sigma$-compact subset is compact.
\end{cor}
{\noindent\bf Proof.} (a) It is easily seen that $\cal{P}$ is an
open filter. If $C_{\infty\cal{P}}(X)$ is a regular ring, then by
Theorem \ref{2.13}, every closed $\cal{P}$-$CG_{\delta}$ is an
open subset. Now let $A$ be a closed $\cal{P}$-$CG_{\delta}$ which
is an open subset. Then by Lemma \ref{2.4}, $A=Z(f)$ for some
$f\in C_{\infty\cal{P}}(X)$. But $X\setminus A=X\setminus
Z(f)=\bigcup_{n=1}^{\infty}\{x: |f(x)|\geq\frac{1}{n}\}$, hence by
\cite[Theorem 3.8.5]{E}, $X\set A$ is a Lindel\"{o}f subset of $X$
, i.e, $A\in\cal{P}$. By Theorem \ref{2.13},
$C_{\infty\cal{P}}(X)$ is a regular ring.

(b) If $A$ is an open locally compact $\sigma$-compact subset,
then $X\setminus A$ is a closed $\cal{P}$-$CG_{\delta}$, where
$\cal{P}$=$\{A\subsetneq X: X\set A$ is compact$\}$ and
$C_{\infty}(X)=C_{\infty\cal{P}}(X)$, so by Theorem \ref{2.13}, we
are done. $\hfill\square$
\\
$$ACKNOWLEDGEMENTS$$
The author would like to thank Prof. Momtahen for his
encouragement and discussion on this paper.


\begin{thebibliography}{99}
\bibitem{AK} S. K. Acharyya and S. K. Ghosh, {\it Functions in $C(X)$ with support lying on
a class of subsets of $X$}, $\textit{Topology procedding}$, $35
(2010), 127$--$148$.
\bibitem{AR} A. R. Aliabad, F. Azarpanah,  M. Namdari,   {\it Rings of continuous functions vanishing at infinity}, $\textit{Comment. Math. Univ. Carolinae}$ 45,3
(2004)$519$--$533$.
\bibitem{ART} A. R.  Aliabad,  F. Azarpanah  and A. Taherifar, \textit{Relative $z$-ideals in commutative rings}, $\textit{comm.
Algebra.}$ $441 (2013) 325$--$341$.
\bibitem{AZ} F. Azarpanah, {\it Intersection of essential ideals
in $C(X)$}, \textit{Proc. Amer. Math. Soc.}, $125(1997),
2149$--$2154$.
\bibitem{A} F. Azarpanah, {\it Essential ideals in $C(X)$},
\textit{Period. Math. Hungar.}, $31(1995), 105$--$112$.
\bibitem{AS} F. Azarpanah, and T. Soundarajan,   {\it When the family
of functions vanishing at infinity is an ideal of $C(X)$},
\textit{Rocky Mountain, journal of mathematics}, $31(4), 2001,
1133$--$1140$.
\bibitem{AT} {F. Azarpanah, and  A. Taherifar}.  Relative $z$-ideals in ${\rm C}(X)$.
\textit{Topology and its applications}. $156 (2009) 1711$--$1717$.
\bibitem{BH}{A. Bella, A. W. Hager, J. Martinez, S. Woodward, H. Zhou, \em Speaker spaces and their
absolutes, I. Top. Appl. $72 (1996)$, $259$--$271$.}
\bibitem{BM}{A. Bella, J. Martinez, S. Woodward, \em Algebra and spaces of
dense constancies, Czechoslovak Math. J. $51 (2001)$,
$449$--$461$.}
\bibitem{E} R. Engelking,  {\it General Topology}, PWN-Polish Sci. Publ  (1977).
\bibitem{G} L. Gilman, and M. jerison, {\it Ring of Continuous Function}, Springer-verlag   (1976).
\bibitem{KM}{M. L. Knox, and W. Wm. McGovern, \em Rigid extensions of $l$-groups of
continuous functions, Czechoslovak Math. J. $58 (2008)$,
$993$--$1014$.}
\bibitem{K} C. W. Kohls,  Ideals in rings of continuous functions, \textit{Fund. Math.} $45 (1957),
28$--$50$.
\bibitem{M} M. Mandelker,   {\it Supports of continuous functions}, \textit{Trans. Amer. Math. Soc.} $156
(1971), 73$--$83$.
\bibitem{MC}{W. Wm. McGovern, \em Clean semiprime $f$-rings with bounded
inversion, Comm. Algebra. $31(2003)$, $3295$--$3304$.}
\end{thebibliography}
\end{document}